\documentclass[12pt]{amsart}

\usepackage{latexsym}
\usepackage{amsfonts}
\usepackage{amssymb}

\def\cal{\mathcal}

\def\bC{{\Bbb C}}
\def\bP{{\Bbb P}}
\def\bF{{\Bbb F}}

\def\cO{{\cal O}}

\def\cC{{\cal C}}

\def\ll{{\frak L}}
\def\ss{\sigma}

\def\dualL{{L\check{\phantom{i}}}}
\def\dualV{{V\check{\phantom{i}}}}

\def\deg{{\operatorname{deg}}}
\def\kod{{\operatorname{kod}}}
\newtheorem{thm}{Theorem}[section]
\newtheorem{lemma}{Lemma}[section]
\newtheorem{proposition}{Proposition}[section]
\newtheorem{remark}{Remark}[section]
\newtheorem{theorem}{Theorem}[section]
\newtheorem{cor}{Corollary}[section]
\newtheorem{example}{Example}[section]
\def\prf{\noindent{\textsc{Proof}}\rm\ }
\def\endprf{\ \hfill $\Box$\medskip}

\newcommand{\grf}[1]{\mbox{$\left \{ #1 \right \}$}}		
\def\unmez{\frac{1}{2}}							
\def\ak{\cO_X(-K_X)}

\begin{document}
\baselineskip=16pt

\parindent=0pt

\title{Classification of Poisson surfaces}

\author[]{Claudio Bartocci and Emanuele Macr\`\i}
{ {\renewcommand{\thefootnote}{\fnsymbol{footnote}}
\footnotetext{\kern-15.3pt AMS Subject Classification: 14J26, 53D17}
}}
\maketitle


\begin{abstract} 
We study complex projective surfaces admitting a Poisson structure; we prove a classification theorem and count how
many independent Poisson structures there are on a given Poisson surface. 
\end{abstract}

\section{Introduction}

The notion of Poisson manifold naturally arises within the framework of analytical
mechanics.  We briefly recall that a Poisson structure on a $C^\infty$ manifold $M$ is given by a bilinear skew-symmetric
bracket $\grf{\cdot,\cdot}$ defined on the sheaf of functions $\cC_M^\infty$, such that
\begin{itemize}
\item $\grf{f,g}= - \grf{g,f}$;
\item $\grf{f,\grf{g,h}}+ \grf{g,\grf{h,f}} + \grf{h,\grf{f,g}}=0$ (Jacobi identity);
\item $\grf{f,gh} = \grf{f,g}h + \grf{f,h}g$.
\end{itemize}
As it was pointed out by Lichnerowicz \cite{L}, the assignment of a Poisson bracket is equivalent to the assignment of a
skew-symmetric bilinear form on the cotangent bundle $T^\ast M$, i.e., a global section $\Pi \in
\Gamma(M,\Lambda^2TM)$,  satisfying the condition
\begin{equation*}[\Pi, \Pi] =0\,,
\end{equation*}
where $[\cdot,\cdot]: \Gamma(M,\Lambda^2 TM)\otimes \Gamma(M,\Lambda^2 TM)\to \Gamma(M,\Lambda^3 TM)$ is the
Schouten-Nijenhuis bracket \cite{V}. Any symplectic manifold $(M,\omega)$ carries a canonical Poisson structure, given by 
$\Pi(\alpha,\beta) = \omega (X_\alpha, X_\beta)$ where $i(X_\alpha)\omega = \alpha$ and $i(X_\beta)\omega = \beta$; the 
condition $[\Pi, \Pi] =0$ is ensured by $d\omega=0$.

The definition of Poisson structure extends in a natural fashion to complex manifolds.  In particular, on a
complex surface $X$ any (holomorphic) global section $\sigma$ of the anticanonical bundle $ \ak = \Lambda^2 T X$ gives
rise to a (holomorphic) Poisson structure, since the condition $[\sigma, \sigma]=0$ is automatically satisfied.   

Complex Poisson surfaces play a major role in the theory of algebraically completely integrable Hamiltonian systems. 
Indeed, as was proved in \cite{B1,B2, Ty}, the choice of a Poisson bivector on a surface $X$
determines natural Poisson structures both on the moduli space of stable sheaves on $X$ and on the Hilbert scheme of points
of $X$. By using this construction, under suitable hypotheses, it is possible to associate an integrable system to a linear
system defined on a Poisson surface, generalizing the results obtained by Beauville for linear systems on K3 surfaces
\cite{B}. Important examples, like the Neumann system, the Hitchin system, etc., can be obtained in this way \cite{Hu,Van}

It is quite immediate to get convinced that a projective Poisson surface can be only an abelian, or a K3 or a ruled
surface. (By ``ruled surface'' we mean any projective surface birationally equivalent to $C\times \bP^1$, where $C$ is a
smooth curve). However, not every ruled surface admits a Poisson structure.  The following classification theorem holds.
\begin{thm}\label{main}
Let $X$ be a minimal ruled surface over the curve $C$ of genus $g$, determined by a normalized rank two vector bundle $V$
over $C$. Let $e = -\deg V$.
\begin{enumerate}
\item If $g=0$, then $X$ is a Poisson surface. 
\item If $g=1$, then
	\begin{itemize}
	\item if $e=-1$, $X$ does not admit any Poisson structure;
	\item if $e \geq 0$, $X$ is a Poisson surface.
	\end{itemize}
\item If $g \geq 2$, then
	\begin{itemize}
	\item if $-g \leq e \leq 2g-3$, $X$ does not admit any Poisson structure;
	\item if $e=2g-2$ and $V$ is indecomposable, $X$ is a Poisson surface;
	\item if $e=2g-2$ and $V$ is decomposable, or $2g-2< e\leq  3g-3$, 
$X$ is a Poisson surface if and only if       
$-K_C - \Lambda^2 V$ is effective;
	\item if $e \geq 3g-2$, $X$ is a Poisson surface.
\end{itemize}
\end{enumerate}
\end{thm}
This theorem can be obtained as a corollary of Sakai's results about the anti-Kodaira dimension of ruled surfaces \cite{S}.
However, Sakai's proof cannot be adapted to answer an important question: how many independent Poisson structures
 are there on a given ruled Poisson surface? In this note we provide a new (and completely elementary) proof of Theorem
\ref{main}; in Corollary \ref{corollario1} and Corollary \ref{corollario2} we compute, whenever it is possible, the
dimension of
$H^0 (X,
\cO_X(-K_X))$.

Related results about $\vert - m K_X\vert$, for an integer $m\geq 1$, can be found in \cite{DZ}.

\section{Classification theorem}

Let $X$ be a smooth projective surface over $\bC$ endowed with a (nontrivial) Poisson structure, namely a nonzero 
section $\ss$ of $\ak$.  Let $D$ be the divisor associated to $\ss$; from the exact sequence 
$$
0 \to \cO_X(K_X) \xrightarrow{\ss} \cO_X \to  \cO_D \to 0\,,$$
it follows that  the Kodaira dimension of $X$
has to be equal either to $0$ or to $-\infty$ (see e.g. \cite[Prop. 2.3]{B1}). In the first case,  
$X$ is a K3 or an Abelian surface, and the canonical bundle is trivial: thus, the Poisson structure is induced by a
(holomorphic) symplectic structure on $X$. If $\kod X = -\infty$,  Enriques' theorem implies that $X$ is a ruled
surface. We notice that, since the section $\ss$ does not vanishes on the open subset $X\backslash D$, the inverse of the
Poisson  bivector is a symplectic form on $X\backslash D$. In other words, $X\backslash D$ is the unique symplectic leaf
of the foliation determined by $\Pi$ \cite{V}. In particular, it follows that Poisson structures on  projective surfaces
have no nontrivial Casimir (holomorphic) functions.   

\begin{example}\rm Let $C$ be a smooth curve. The cotangent bundle $T^\ast C$ carries a canonical symplectic form
$\Omega=d\theta$, where $\theta$ is the Liouville one-form. Denoting by $Q$ the total space of $T^\ast C$, it follows that
$Q$ is a non-compact symplectic surface. It is easy to show that $Q$ can be embedded as an open set into the ruled surface
$X= \bP(\cO_C \oplus \cO_C(K_C))$. $X$ is a Poisson surface and $Q$ a symplectic leaf.
\end{example}
It is obvious, however, that not every ruled surface carries a Poisson structure. For instance, the anticanonical bundle
of the surface
$X= C\times \bP^1$ has no nonzero sections unless $g(C)\leq 1$. 

In order to classify the ruled surfaces admitting Poisson
structures, we have first to understand what happens when a Poisson surface $X$ is blown-up at a point $p$.  

\begin{lemma} Let $\rho: \tilde X \to X$ the blow-up at the point $p\in X$. Then
\begin{equation}\label{ineq}
h^0 (X, \cO_X (-K_X)) -1 \leq h^0 (\tilde X , \cO_{\tilde X} (-K_{\tilde X})) \leq h^0 (X, \cO_{X} (-K_X))\,,
\end{equation}
and $h^0 (\tilde X , \cO_{\tilde X} (-K_{\tilde X})) = h^0 (X, \cO_{X} (-K_X))$ if and only if $p$ is a base point of 
$\vert -K_X \vert$. 
\end{lemma} 
\prf
Let $E$ be the exceptional divisor. Since $\rho$ is an isomorphism on $\tilde X\backslash E$, 
any Poisson bracket on $\tilde X$ induces a Poisson bracket on $X$ by Hartogs' theorem; this proves the right inequality 
in (\ref{ineq}). Since $-K_{\tilde X} \cdot E >0$, any section  $\ss\in H^0 (X, \cO_{X} (-K_X))$ coming from a section  
$\tilde \ss \in H^0 (\tilde X , \cO_{\tilde X} (-K_{\tilde X}))$ passes through the point $p$. Conversely, if $\ss$ passes
through the point $p$, then $-K_{\tilde X} \sim \rho^\ast D$, where $D$ is the divisor associated to $\ss$, and therefore 
$-K_{\tilde X}$ is effective. 
\endprf 

We can now restrict our attention to minimal ruled surfaces. We shall freely use the results and  
notations in \cite{H}, Chap. V, \S 2. 

If $q(X)=0$, then $X$ is $\bP^2$ or the rational ruled surface $\bF_n$, with $n\neq 1$. 
In both cases, a straighforward computation shows that $h^0(X, \cO_{X} (-K_{X})) \geq 9$.
\begin{proposition} Any minimal ruled surface $X$ with $q(X)=0$ is a  Poisson surface.
\endprf
\end{proposition}

When $q(X) \geq 1 $, then $\pi: X\to C$ is a geometrically ruled surface, with $g(C)=q(X)$. 
We can assume $X = \bP(V)$, where $V$ is rank $2$ vector bundle on $C$ such that $H^0(C,V) > 0$ and $H^0(C, V\otimes M)=0$
for every line bundle $M$ of negative degree. We shall say that such an $V$ is {\it normalized}. Under this hypothesis, 
there exists a section $\tau : C\to X$ such that 
$\tau^2= \deg V =: - e$; we have $\cO_X(\tau) \cong \cO_X(1)$ and so $V = R^0\pi_\ast \cO_X(\tau)$. 
There are some restrictions on the possible values of the invariant $e$ \cite{N,H}: 
\begin{enumerate}
\item if $V \cong \cO_C \oplus L$, then $e\geq 0$;
\item if $V$ is indecomposable, then $- g \leq e \leq 2g-2$. 
\end{enumerate}
Moreover, all these values are admissible.

Any normalized vector bundle $V$ over the curve $C$ fits into an exact sequence
\begin{equation}\label{seq1}
0 \to \cO_C \to V \to L \to 0,
\end{equation}
where $L$ is a line bundle over $C$; we have   $L\cong\Lambda^2 V$. Let $\ll$ be the divisor on $C$ corresponding 
to $L$; it is easy to show that
\begin{equation}\label{21}
-K_X \sim 2\tau  + \pi^* (-K_C - \ll)\,.
\end{equation}
By the projection formula, we obtain
\begin{equation}\label{iso}
H^0(X, \cO(-K_X) )\cong H^0(C,S^2(V)\otimes \dualL(-K_C))\,.
\end{equation}

\begin{lemma}\label{lemmaseq}
There is an exact sequence of vector bundles over $C$:
\begin{equation}\label{seq2}
0 \to V \otimes \dualL(-K_C) \to S^2 (V) \otimes \dualL(-K_C) \to L(-K_C) \to 0\,.
\end{equation}
\end{lemma}
\prf
Let us consider the exact sequence
$$0 \to \cO_X (1) \to \cO_X (2) \to \cO_X(2) \otimes \cO_{\tau} \to 0\,;$$
since $R^1 \pi_* (\cO_X (1)) =0$ (see \cite{H}, Chap. V, Lemma 2.4), we get the exact sequence
$$0 \to V \to S^2 (V) \to L \otimes L \to 0\,.$$
The result follows by tensoring this sequence by the line bundle $\dualL(-K_C)$.
\endprf
\begin{remark}\rm By using the exact sequence (\ref{seq2}) it is an easy exercise to compute the number of independent
Poisson structures on the surfaces $\bF_n$ (in this case, the invariant $e$ coincides with $n$):
\begin{itemize}
\item if $X=\bF_0 \cong \bP^1\times \bP^1$, then $h^0 (X , \cO (-K_{X})) = 9$; 
\item if $X=\bF_2$, then $h^0 (X , \cO (-K_{X})) = 9$;
\item if $X=\bF_n$, $n\geq 3$, then $h^0 (X , \cO (-K_{X})) = n+6$.
\end{itemize}\endprf
\end{remark}

We recall a useful criterion of ampleness that we shall exploit in the proof of Proposition \ref{proplemma}.

\begin{proposition}\label{2.3}
\begin{enumerate}
\item If $e \geq 0$, a divisor $D \equiv a\tau + bf$ is ample if and only if $a > 0$ and $b > ae$.
\item If $e < 0$, a divisor $D \equiv a\tau + bf$ is ample if and only if  $a > 0$ and $b> \unmez ae$.
\end{enumerate}\endprf
\end{proposition}

\begin{theorem}\label{theor1} Let $X$ be a minimal ruled surface with $q(X)=1$. 
\begin{enumerate}
\item If $e=-1$, $X$ does not admit any Poisson structure;
	\item if $e \geq 0$, $X$ is a Poisson surface.
\end{enumerate}
\end{theorem}
\prf

If $e=0$, we have to distinguish 3 cases: $V$ indecomposable, $V= \cO_C \oplus \cO_C$ and 
$V= \cO_C \oplus L$ with $L\ncong \cO_C$ (recall that $V$ is  normalized). If $V$ is indecomposable, then, as shown
in Theorem V.2.15 of \cite{H}, it is uniquely determined by the exact sequence
\begin{equation}\label{seq3}
0 \to \cO_C \to V \to \cO_C \to 0\,.
\end{equation}
Since the sequence (\ref{seq3}) does not split, we have $h^0(C,V)= 1$. 
The long exact cohomology sequence associated to the
exact sequence (\ref{seq2}) is
\begin{equation}\label{g=0e=0} 
0 \to H^0(C,V) \to H^0(C, S^2(V)) \to \bC \to \bC \to  H^1(C, S^2(V)) \to \bC \to 0\,.
\end{equation}
Now, $S^2(V)$ is indecomposable (see \cite[Th. 9]{A}), so that
one has  $h^0(C, S^2(V)) = h^0(X, \cO(-K_X))= 1$. The case $V= \cO_C \oplus \cO_C$, corresponding to $X = C\times
\bP^1$, is trivial; we have
$$H^0 (X, \cO(-K_X)) \cong H^0 (C, \cO_C) \otimes H^0 ( \bP^1 , \cO_{{\bP}^1} (2))\,,$$
hence, $h^0 (X, \cO(-K_X)) = 3$. If $V= \cO_C \oplus L$ with $L\ncong \cO_C$, then $h^0(C, L)= h^1(C, \dualL)=0$. From
the exact sequence
$$0 \to \dualL \to V \otimes \dualL \to \cO_C \to 0\,,$$
we get $H^0 (C, V \otimes \dualL) \cong H^0 (C, \cO_C) \cong \bC$. By using again the sequence (\ref{seq2}), it follows
$h^0 (X, \cO(-K_X)) = 1$.

If $e> 0$, then $V$ is decomposable: $V \cong \cO_C \oplus L$, with $\deg L <0$. By reasoning as in the previous
case, we get  $h^0 (X, \cO(-K_X)) = e+1$.

Finally, if $e=-1$, then $V$ is uniquely determined by the exact sequence
\begin{equation}\label{seq4}
0 \to \cO_C \to V \to \cO_C (P) \to 0\,,
\end{equation}
where $P$ is a point of $C$. Dualizing the sequence (\ref{seq4}), we get at once that $\dualV \cong V\otimes \cO_C(-P)$.
Now, $V$ is indecomposable, of rank
$2$ and degree
$1$; so, by
\cite[Lemma 22]{A}, one has
$V\otimes \dualV \cong \oplus_{i=0}^3 \Xi_i$, where the $\Xi_i$ are the line bundles on $C$ of order diving $2$ (in
particular, we set $\Xi_0 \cong \cO_C$). An easy computation shows that $S^2 (V)\otimes \cO_C(-P)\cong \oplus_{i=1}^3
\Xi_i$. One concludes that  
$h^0 (X,
\cO (-K_X))= h^0(C, S^2 (V)\otimes \cO_C(-P)) = 0$.
\endprf

\begin{cor}\label{corollario1} Let $X$ be a minimal ruled surface with $q(X)=1$.
	\begin{enumerate}
	\item If $e=0$ and $V = \cO_C\oplus \cO_C$ (hence $X \cong \bP^1 \times C$), then $h^0 (X, \cO_X(-K_X))=3$; 
	\item if $e=0$ and $V$ is indecomposable, then  $h^0 (X, \cO_X(-K_X)) =1$;
	\item if $e=0$ and $V = \cO_C \oplus L$,with $L \neq \cO_C$, then $h^0 (X, \cO_X(-K_X)) =1$;
	\item if $e \geq 1$, then $h^0 (X, \cO_X(-K_X)) = e+1$.
	\end{enumerate}
\endprf
\end{cor}

It follows from the proof of Theorem \ref{theor1} that, for a ruled surface over an elliptic curve, it may happen that the
divisor $-\ll$ corresponding to $-\Lambda^2 V$ is not effective while $-K_X \sim 2\tau + \pi^\ast(-\ll)$ is
effective. This is not the case when $g(C) >1$, as we shall prove in the following Proposition, which can be rephrased as
follows:  the divisor $-K_X\sim 2\tau  +\pi^* (-K_C -\ll)$ is effective if and only if $- K_C -\ll$ is effective, where
$\ll$ is the divisor corresponding to the line bundle $\Lambda^2 V$.

\begin{proposition}\label{proplemma} Let $\pi :X \to C$ be a minimal ruled surface with $q(X)=g(C) >1$. Then,
$h^0(X,\cO_X(-K_X)) = h^0(C,
\dualL(-K_C))$.
\end{proposition}
\prf
From the exact sequence (\ref{seq2}) we obtain:  
\begin{equation}
0 \to H^0(C,V \otimes \dualL(-K_C) )\to H^0(C, S^2 (V) \otimes \dualL(-K_C)) \to H^0(C, L(-K_C))\,.
\end{equation}
Now, $\deg L(-K_C) = - e -2g +2$; since $-e \leq g$, one has $\deg L(-K_C) <0$ in all cases except when $g=2$, $e=-2$.
So, for $g\geq 3$ and for $g=2$, $e\neq -2$, we get
$$h^0(C,V \otimes \dualL(-K_C) ) = h^0(C, S^2 (V) \otimes \dualL(-K_C))\,.$$
The exact sequence (\ref{seq1}) implies $h^0(C,V \otimes \dualL(-K_C) ) = h^0(C,\dualL(-K_C) )$; thus, by (\ref{iso}) we
get 
$h^0(X,\cO_X(-K_X)) = h^0(C, \dualL(-K_C)$. 
To deal with the missing case $g=2$, $e=-2$, we use the criterion in Proposition \ref{2.3}: the divisor $\tau$ is ample.
But,
$ -K_X\cdot \tau =0$, so $-K_X$ is not effective. This ends the proof.
\endprf

We can make the previous statement somewhat more precise. By noticing that if $e=2g-2$ and $V$ is
indecomposable, then $h^1(C,\dualL) \neq 0$, hence $h^0(C, L(K_C)) \neq 0$, it is indeed easy to prove the following
result.  

\begin{cor}\label{corollario2} Let $X$ be a minimal ruled surface with $q(X)=g \geq 2$. 
	\begin{enumerate}
	\item If $-g \leq e \leq 2g-3$, $X$ does not admit any Poisson structure;
	\item if $e=2g-2$ and $V$ is indecomposable, $X$ is a Poisson surface;
	\item if $e=2g-2$ and $V$ is decomposable, or $2g-2< e\leq  3g-3$, 
$X$ is a Poisson surface if and only if       
$-K_C - \Lambda^2 V$ is effective;
	\item if $e \geq 3g-2$, $X$ is a Poisson surface.
	\end{enumerate}\endprf
\end{cor}
We recall that, whenever $2g-2< e\leq  3g-3$, the vector bundle $V$ is decomposable.

\vskip20pt
\paragraph{\bf Acknowledgements.} 
C.B. acknowledges the financial support of the MIUR and the University
of Genova through the national research project ``Geometria dei sistemi integrabili''. The authors thank the referee
for useful remarks, and Arnauld Beauville for having drawn Sakai's paper to their attention after the first draft of this
paper was completed.

\vfill

\hrule
\smallskip

\bigskip\baselineskip=2pt
\footnotesize{\obeylines
\noindent Claudio Bartocci
\noindent Dipartimento di Matematica
\noindent Universit\`a degli Studi di Genova
\noindent Via Dodecaneso 35
\noindent 16146 Genova, ITALY
\noindent E-mail: bartocci@dima.unige.it
\noindent Emanuele Macr\`\i
\noindent Scuola Internazionale di Studi Superiori Avanzati (SISSA)
\noindent Via Beirut 2-4
\noindent 34014 Trieste, ITALY
\noindent E-mail: macri@sissa.it}

\end{document}